\documentclass{amsart}
\usepackage{amsmath, amsthm, amscd, amsfonts, amssymb, graphicx, color}
\usepackage[bookmarksnumbered, colorlinks, plainpages]{hyperref}

\newtheorem{theorem}{Theorem}[section]
\newtheorem{cor}[theorem]{Corollary}

\newtheorem{pro}[theorem]{Proposition}
\theoremstyle{definition}

\theoremstyle{Example}

\numberwithin{equation}{section}
\begin{document}
\title [ Euclidean Geometry and  Curves ]
{Euclidean Geometry and Elliptic Curves }
\author[Farzali   Izadi]{Farzali Izadi}
\address{Department of Mathematics, University Of Toronto, Toronto ON Canada}
\email{\href{f.izadi@utoronto.ca}{f.izadi@utoronto.ca} }
\date{October, 13, 2019}

%\date{January 5, 2011 and, in revised form, June 5, 2011.}

\maketitle
\begin{abstract}
In this paper, we demonstrate the intimate relationships among some geometric figures and the families of elliptic curves with positive ranks. These geometric figures include \textit{\textbf{Heron triangles}}, \textit{\textbf{Brahmagupta quadrilaterals}} and \textit{\textbf{Bicentric quadrilaterals}}.
Firstly, we investigate the important properties of these figures and then utilizing these properties, we show that how to construct various families of elliptic curves with different positive ranks having different torsion subgroups.
\end{abstract}
\section{Introduction}
\label{intro} In geometry, {\textbf{Heron$^{,}$s formula} (sometimes called \textbf{Hero$^{,}$s formula}), named after hero of Alexandria, gives the area of a triangle by requiring no arbitrary choice of side as base or vertex as origin, contrary to other formulas for the area of a triangle, such as half the base times the height or half of the norm of a cross product of two sides.\\
Precisely, Heron$^{,}$s formula states that the area of a triangle whose sides have lengths $a$, $b$, and $c$ is
 $$A=\sqrt{p(p-a)(p-b)(p-c)},$$
where $p$ is the semiperimeter of the triangle; that is,
$$p=\frac{a+b+c}{2}.$$

The Indian mathematician Brahmagupta, 598-668 A.D.,
showed that for a triangle with integral sides ${a, b, c}$
and integral area $A$ there
are positive integers $k, m, n$, with $k^2 < mn$ , such that
\begin{align}
\begin{cases}
a&=n(k^2 + m^2),\\
b&=m(k^2 + n^2),\\
c&=(m + n)(mn − k^2),\\
A&=kmn(m + n)(mn − k^2).
\end{cases}
\end{align}
We observe that the value of $A$ is a consequence of the Heron formula. A classical
reference for Heron triangles is the second volume of the History of the Theory of
Numbers by L. E. Dickson \cite{Dic}.\\

Heron$^{,}$s formula is a special case of \textbf{Brahmmagupta formula} for the area of a \textbf{cyclic quadrilateral} (i.e., a simple quadrilateral that is inscribed in a circle) whose sides, diagonals and area are with sides of length $a$, $b$, $c$, and $d$ as
$$A=\sqrt{(p-a)(p-b)(p-c)(p-d)},$$
where $p$ is the semiperimeter
$$s=\frac{a+b+c+d}{2}.$$
Similarly, a \textbf{bicentric quadrilateral} is a \textbf{convex quadrilateral} that has both a circumcircle passing through the four vertices and incircle having the four sides as tangents. In what follows, we consider a bicentric quadrilateral with rational sides, and discuss the problem of finding such quadrilateral where the ratio of the rdii of the circumcircle and incircle is rational.\\
The radii of these circles are denoted by $R$ and $r$ respectively. Examples of bicentric quadrilateral are squares, right kits, and isosceles tangential trapezoids. By Pitot$^{,}s$ theorem we have $s=a+c=b+d$ where $s$ is the semiperimeter. In addition we have
$$R=\frac{1}{4}\sqrt{\frac{(ab+cd)(ac+bd)(ad+bc)}{abcd}},$$
and $r=\frac{K}{s}$ where $K=\sqrt{abcd}$ is the area of the quadrilateral.\\
In this review paper, we show that many authors utilize these formulas  form all three figures to construct different elliptic curves with positive ranks. The main goal of the paper is that how the Euclidean geometry intimately related to the theory of elliptic curves.
\\
%----------------------------------------------------------------------------------------------------------------
The use of Heron$^{,}$s formula in umber theory was not new. However the author together with some other people provided with a series of papers utilizing Heron$^{,}$s formula, Brahmagupta quadrilaterals as well as Bicentric quadrilaterals in an attempt to construct new families of elliptic curves with positive ranks. In this paper, we try to demonstrate the main results of these series of papers and leave the details for the corresponding original papers. To make it more accessible for the readers we devote separate sections to separate papers. Each section comes with the title of the corresponding original paper.

\section{A Family of Elliptic Curves With Rank $\geq5$}
In this paper, \cite{IN2} the authors construct a family of elliptic curves with rank $\geq 5$. To do this, they use the  Heron formula for a triple $(A^2, B^2, C^2)$ which are not necessarily the three sides of a triangle. It turns out that as parameters of a family of elliptic curves, these three positive integers $A$, $B$, and $C$, along with the extra parameter $D$ satisfy the quartic Diophantine equation $A^4+D^4=2(B^4+D^4)$. For the solutions of this diophantine equation see \cite{IK}.

Take  $(A^2, B^2, C^2)$, where   $(A, B, C)$  are as in diophantine equation mentioned above. Since the triple $(A, B, C)$ is arising from the Diophantine equation, there is no guarantee to have a real triangle. Now by taking

       $a=A^2$, $b=B^2$, and $c=C^2$, one finds that
\begin{equation}\label{Eq22}
S=\sqrt{\frac{(A^2+B^2+C^2)(A^2+B^2-C^2)(A^2+C^2-B^2)(B^2+C^2-A^2)}{16}}.
\end{equation}
Expanding \eqref{Eq22}, one gets
\begin{equation*}
S^2=-\left(\frac{A^8+B^8+C^8-2A^4B^4-2A^4C^4-2B^4C^4}{16}\right),
\end{equation*}
equivalently,
\begin{equation*}
16S^2=2A^4B^4+2A^4C^4+2B^4C^4-A^8-B^8-C^8,
\end{equation*}
or
\begin{equation*}
\left(\frac{A^4+B^4-C^4}{2}\right)^2+4S^2=A^4B^4.
\end{equation*}

\begin{equation*}
\left(\frac{A^4+B^4-C^4}{2}\right)^2+4S^2=A^4B^4.
\end{equation*}
Multiplying both sides by $A^2B^2$ yields
\begin{equation*}
A^2B^2\left(\frac{A^4+B^4-C^4}{2}\right)^2+4A^2B^2S^2=A^6B^6.
\end{equation*}
Taking $y=AB\left(\frac{A^4+B^4-C^4}{2}\right)$ and $x=A^2B^2$, the following  family of elliptic curves is obtained:
\begin{equation}\label{Eq1.2}
E:\ y^2=x^3-4S^2x.
\end{equation}

Since the roles of $A$, $B$, and $C$ are symmetric in the Heron formula, the following points lie on
the family \eqref{Eq1.2} also.
\begin{equation*}
\begin{array}{l}
P_1=\left(A^2B^2,\ \frac{AB(A^4+B^4-C^4)}{2}\right),\\
P_2=\left(A^2C^2,\ \frac{AC(A^4+C^4-B^4)}{2}\right),\\
P_3=\left(B^2C^2,\ \frac{BC(B^4+C^4-A^4)}{2}\right),
\end{array}
\end{equation*}
Next it is easy to impose two more  points on the curve \eqref{Eq1.2} with $x-$coordiates as $B^2D^2$ and $C^2D^2$. Substituting $x=B^2D^2$ in \eqref{Eq1.2} yields
\begin{equation*}
y^2=B^2D^2\left(\frac{4B^4D^4+A^8+B^8+C^8-2A^4B^4-2A^4C^4-2B^4C^4}{4}\right),
\end{equation*}
or
\begin{equation}\label{Eq1.3}
y^2=B^2D^2\left(\frac{A^4(A^4-2B^4-2C^4)+B^8+C^8-2B^4C^4+4B^4D^4}{4}\right).
\end{equation}
Let
\begin{equation}\label{Eq1.4}
A^4-2B^4-2C^4=-D^4,
\end{equation}
then by substituting $A^4=2B^4+2C^4-D^4$ in \eqref{Eq1.3}, we get
\begin{equation*}
y^2=B^2D^2\left(\frac{B^2+D^2-C^2}{2}\right)^2.
\end{equation*}

Thus, the point $P_4=\left(B^2D^2,\ \frac{BD(B^4+D^4-C^4)}{2}\right)$ is a new  point on \eqref{Eq1.2}.
Similarly, one can easily check that the point $P_5=\left(C^2D^2,\ \frac{CD(C^4+D^4-B^4)}{2}\right)$ lies  on \eqref{Eq1.2} as well.\\

It is clear that the existence  of these extra points on the family depends exactly on the existence of the solutions of the Diophantine equation $A^4+D^4=2(B^4+C^4)$.\\

In the final stage, the authors proved that the torsion subgroup of the \eqref{Eq1.2} is $\Bbb{Z}/2\Bbb{Z}$.

\section{Rank of elliptic curves associated to Brahmagupta quadrilaterals}

In this paper, \cite{IKZ1} the authors construct a family of elliptic curves with six parameters, arising from a system of Diophantine equations, whose rank is at least five. To do so, they use the Brahmagupta formula for the area of cyclic quadrilaterals \((p^3,q^3,r^3,s^3)\) not necessarily representing genuine geometric objects. It turns out that, as parameters of the curves, the integers \(p,q,r,s\), along with the extra integers \(u\), \(v\) satisfy the system of diophantine equations \(u^6+v^6+p^6+q^6=2(r^6+s^6)\), \(uv=pq\), (see \cite{IKZ2}), and they then utilize a subset of solutions of the system by the solutions of a specific positive rank elliptic curve.

In this work, the authors use the Brahmagupta formula for the area of cyclic quadrilaterals to find infinitely many elliptic curves with high rank. In effect, the elliptic curve is
$$
y^2=x^3-3(pqrs)^2x+2(pqrs)^3+\frac{1}{4}(p^6+s^6-q^6-r^6)^2-(p^3s^3+q^3r^3)^2,
$$
denoted by $E_{u,v,p,q,r,s}$, over
\begin{align}
& u^6+v^6+p^6+q^6=2(r^6+s^6) \label{eq:1} \\
C: \quad & \nonumber \\
& uv=pq \label{eq:2}
\end{align}
and prove that the group of the rational maps $C\rightarrow E_{u,v,p,q,r,s}$, that commute with the projection $E_{u,v,p,q,r,s}\rightarrow C$, has rank at least five. This is done by exhibiting five explicit sections $P_1$, $P_2$, $P_3$, $P_4$, $P_5$, and showing that these are linearly independent. They use the fact that $C(\Bbb Q)$, the set of rational solutions on $C$, is infinite in order to deduce that infinitely many specializations of $E_{u,v,p,q,r,s}$ have ranks at least five over the rationales. This last assertion is done by using $$y^2-28xy-560y=x^3-20x^2-400x+8000,$$
an elliptic curve of positive rank lying on $C$ found in \cite{IKZ2}. In fact, the following theorem is proved.

\begin{theorem}
There are infinitely many elliptic curves $E_{u,v,p,q,r,s}$ over $C$ of rank at least five, parameterized by an elliptic curve of rank at least three over $\Bbb Q(p,q,r,s)$.
\end{theorem}
\subsection{The construction of $E_{u,v,p,q,r,s}$}
As is mentioned above, in this work, the authors deal with the elliptic curves related to (positive) integer solutions on $C$. The Brahmagupta formula for the area of a cyclic quadrilateral in terms of its side lengths states that for a quadrilateral with sides $a,b,c,d$, the area of the quadrilateral is given by
\begin{equation}\label{hq1}
S=\sqrt{(P-a)(P-b)(P-c)(P-d)},
\end{equation}
in which $P=(a+b+c+d)/2$ is the semi-perimeter. Equivalently,
\begin{equation}\label{hq2}
16S^2=(b+c+d-a)(a+c+d-b)(a+b+d-c)(a+b+c-d),
\end{equation}

which can be written in the form
\begin{equation}\label{hq3}
\frac{1}{4}\left(b^2+c^2-a^2-d^2\right)^2=(ad+bc)^2-4S^2.
\end{equation}
Now, take $(a,b,c,d)=(p^3,q^3,r^3,s^3)$. By substitution, \eqref{hq3} turns into
\begin{equation}\label{hq4}
\frac{1}{4}\left(q^6+r^6-p^6-s^6\right)^2=(p^3s^3+q^3r^3)^2-4S^2.
\end{equation}
Expanding and sorting the right hand side of \eqref{hq4}, leads to
$$\frac{1}{4}\left(q^6+r^6-p^6-s^6\right)^2=(p^2s^2)^3+(q^2r^2)^3+2(pqrs)^3-4S^2,$$

or
$$\left(\frac{q^6+r^6-p^6-s^6}{2}\right)^2=(p^2s^2+q^2r^2)^3-3p^2q^2r^2s^2(p^2s^2+q^2r^2)+2(pqrs)^3-4S^2.$$
Setting $x=p^2s^2+q^2r^2$, $y=(q^6+r^6-p^6-s^6)/2$, and using \eqref{hq4}, one can define the following elliptic curve
\begin{equation}\label{curve0}
y^2=x^3-3(pqrs)^2x+2(pqrs)^3+\frac{1}{4}(p^6+s^6-q^6-r^6)^2-(p^3s^3+q^3r^3)^2,
\end{equation}
denoted by $E_{p,q,r,s}$, over $\Bbb Q(p,q,r,s)$.

By symmetric roles of $p, q, r, s$ in the formulas, the elliptic curve has the following non-obvious points:
\begin{align*}
P_1(p,q,r,s)&=\left(p^2q^2+r^2s^2,\frac{p^6+q^6-r^6-s^6}{2}\right), \\
P_2(p,q,r,s)&=\left(p^2r^2+q^2s^2,\frac{p^6-q^6+r^6-s^6}{2}\right), \\
P_3(p,q,r,s)&=\left(p^2s^2+q^2r^2,\frac{p^6-q^6-r^6+s^6}{2}\right).
\end{align*}
By the specialization theorem \cite{S2}, they proved that the family has rank at least 3. Finally utilizing the
system of equations
$$u^6+v^6+p^6+q^6=2(r^6+s^6), uv=pq.$$
discussed in \cite{IKZ2} the authors showed that the following points also lie on the above family of elliptic curves.
\begin{align}
P_4(u,v,p,q,r,s)&=\left(u^2r^2+v^2s^2,\frac{u^6-v^6+r^6-s^6}{2}\right), \\
P_5(u,v,p,q,r,s)&=\left(u^2s^2+v^2r^2,\frac{u^6-v^6-r^6+s^6}{2}\right),
\end{align}

\section{On Elliptic Curves Via Heron Triangles and Diophantine Triples}
In this paper F. Izadi et al. \cite{IK}, constructed families of elliptic curves arising from Heron triangles and Diophantine triples with the Mordell-Weil torsion subgroup of ${Z/2Z}\times {Z/2Z}$. These families have ranks at least 2 and 3, respectively, and contain particular examples with rank equal to 7.\\
The following are the main results:

Consider the elliptic curve $E_{k}: y^{2}=(x+a(k)b(k))(x+b(k)c(k))(x+a(k)c(k))$ associated to the triples:
\begin{align}
\begin{cases}
a(k)&=10k^{2}-8k+8, \\
  b(k)&=k(k^{2}-4k+20),\\
 c(k)&=(k+2)(k^{2}-4),
\end{cases}
\end{align}

%-----------------------------------------------------------------------------------------
\begin{theorem} Let $a(k)$, $b(k)$ and $c(k)$ be defined as (4.1) functions, where $k$ is an arbitrary rational number different from $0$, $-2$ and $2$. Then the elliptic curve
$$E: y^{2}=(x+a(k)b(k))(x+b(k)c(k))(x+a(k)c(k))$$
defined over $\mathbb{Q}(k)$ has torsion subgroup ${\mathbb{Z}/2\mathbb{Z}}\times {\mathbb{Z}/2\mathbb{Z}}$.
\end{theorem}
%----------------------------------------------------------------------------
\begin{theorem}
With the terminology in theorem (4.1), rank $E(\mathbb{Q}(k))\geqslant 2$.
\end{theorem}

\begin{pro}
For each $2\leqslant r\leqslant 7$, there exists some $k$ such that the elliptic curve $E_{k}$ defined in the theorem (4.1) has torsion subgroup $\mathbb{Z}/2\mathbb{Z}\times \mathbb{Z}/2\mathbb{Z}$ and rank $r$.
\end{pro}

\begin{theorem}
There exists a subfamily of $E_{k}$ of rank $\geqslant 3$ over $\mathbb{Q}(m)$.
\end{theorem}

The last result related to the diophantine triples.
Let the diophantine triple $(a,b,c)$ has the property $D(n)$ (\cite{DP2}) if for any non zero integer $n$, there exist rationales $r$, $s$ and $t$ such that $$ab+n=r^{2},\quad \quad \quad  ac+n=s^{2}, \quad \quad \quad  bc+n=t^{2}.$$
%------------------------------------------------------------------------------------------------
\begin{theorem}
Consider $(a,b,c)=(k-1,k+1,4k)$ with property $D(1)$. Then there exists a subfamily of $C:Y^{2}=(aX+1)(bX+1)(cX+1)$ over $\mathbb{Q}$ with rank $\geqslant 2$.
\end{theorem}

\section{Elliptic Curves Coming From Heron Triangles}
The paper by Dujella and peral \cite{DP1} adopted the method  from \cite{IK} to construct families of elliptic curves with ranks at least $2$, $3$, $4$, and $5$. For the subfamily with rank $\geqslant 5$, they show that its generic rank is exactly equal to 5 and found free generators of the corresponding group. By specialization,
examples of elliptic curves over $Q$ with rank equal to 9 and 10 are also obtained.
 \\
\section{Heron Triangles Via Elliptic Curves}
   Goin and Maddox \cite{GM} generalized some of Koblit$^{,}$s notions \cite{Koblitz}, [Ch. 1, \S 2, ex. 3] by exploring the correspondence between positive integers $n$ associated with arbitrary triangles (with rational side lengths) which have area $n$ and the family of elliptic curves

\begin{align}
 E_{\tau}^{(n)}:y^{2}=x(x-n\tau)(x+n\tau ^{-1}).
\end{align}

for nonzero rational $\tau$. Congruent number elliptic curves are obtained by taking $\tau=1$. In fact, they broadened the study from rational right triangles to simply arbitrary rational triangles known as the Heron triangles.

The main result of the Goin-Maddox paper is the following.

\begin{theorem} A positive integer $n$ can be expressed as the area of a triangle with rational sides if and only if for some nonzero rational number $\tau$ the elliptic curve
 \begin{align}
 E_{\tau}^{(n)}:y^{2}=x(x-n\tau)(x+n\tau ^{-1})
\end{align}
has a rational point which is not of order 2. Moreover, $n$ is a congruent number if and only if we can choose $\tau=1$.
\end{theorem}

The more general elliptic curve $E_{\tau}^{(n)}$ constructed above has all of its 2-torsion rational. Explicitly, when
$\tau=4n/((a + b)^{2}-c^{2})$, then
 \begin{align}
 E_{\tau}^{(n)}[2]=\begin{Bmatrix}
 (0,0), & (\frac{c^{2}-(a + b)^{2}}{4},0), & (\frac{c^{2}-(a-b)^{2}}{4},0), & \mathcal{O}
\end{Bmatrix}.
\end{align}\\

They showed that the torsion subgroup of $E_{\tau}^{(n)}$ is either
 \begin{align}
  \mathbb{Z}/2\mathbb{Z}\times \mathbb{Z}/2\mathbb{Z} \quad or \quad \mathbb{Z}/2\mathbb{Z}\times \mathbb{Z}/4\mathbb{Z},
\end{align}

with the latter case corresponding to isosceles triangles.

\section{Heron Quadrilaterals Via Elliptic Curves}
In other direction, Izadi, Khoshnam and Moody, \cite{IKM}  established a correspondence between cyclic quadrilaterals with rational side lengths and area $n$,
i.e., the correspondence between the Heron quadrilaterals and the family of elliptic curves of the form
 \begin{align}
y^2=x^3+\alpha x^2-n^2 x
\end{align}
This correspondence generalizes the notions of Goins-Maddox \cite{GM} who established a similar connection between Heron triangles and elliptic curves.
The authors further study this family of elliptic curves, looking at their torsion groups and ranks. They also explore their connection with congruent numbers, which are the $\alpha=0$ case. Congruent numbers are positive integers which are the area of a right triangle with rational side lengths.
We give explicit formulas which show how to construct the elliptic curve and some non-trivial points on the curve given the side lengths and area of the quadrilateral.  If we set one of the side lengths to zero, then the formulas collapse to exactly those of Maddox-Goins. They also show the other direction of the correspondence, that is, how to find a cyclic quadrilateral which corresponds to a given elliptic curve in our family. The pair $(\alpha, n)$ is called a \emph{generalized congruent number pair} if the elliptic curve $y^2=x^3+\alpha x^2-n^2 x$ has a point of infinite order. Similarly one can call the curve a generalized congruent number elliptic curve.  The generalized congruent number curves with $\alpha=0$ are precisely the congruent number curves.  Stated in this way, the results relate generalized congruent number pairs with cyclic quadrilaterals with area $n$.
Some important results are as follows:
\begin{theorem}
\label{main}
For every cyclic quadrilateral with rational side lengths and area $n$, there is an elliptic curve
$$E_{\alpha,-n^2}:y^2=x^3+\alpha x^2-n^2x$$
with 2 rational points, neither of which has order 2.  Conversely, given an elliptic curve $E_{\alpha,-n^2}$ with positive rank, then there is a cyclic quadrilateral with area $n$ whose side lengths are rational (under the correspondence given above).
\end{theorem}

\begin{theorem}
\label{rank0}
If the curve $E_{\alpha,-n^2}$ arising from a cyclic quadrilateral has rank 0 then the associated quadrilateral is either a square, or an isosceles trapezoid with three sides equal ($a=b=d$) such that $(a+c)(3a-c)$ is a square.   The torsion group is $\mathbb{Z}/6\mathbb{Z}$ for this rank 0 case.
\end{theorem}

\begin{theorem}
If the quadrilateral is a square, then the rank is  0, with torsion group $\mathbb{Z}/6\mathbb{Z}$.
\end{theorem}

The rank need not be 0 for isosceles trapezoids with three sides equal.  Take, for example the quadrilateral with side lengths $(13,13,23,13)$, which yields the curve $E_{-11,-216^2}$.  This curve has rank 1, with generating point $(-196,1092)$ and has torsion group $\mathbb{Z}/6\mathbb{Z}$. One can have rank 0 curves $E_{\alpha,-n^2}$ with torsion group $\mathbb{Z}/2\mathbb{Z} \times \mathbb{Z}/4\mathbb{Z}$.  Take for example $\alpha=7, n=12$.  The previous results prove that such curves do not correspond with a cyclic quadrilateral. By allowing the quadrilaterals with $d=0$, then it can be shown that these curves come from quadrilaterals with $d=0$, i.e. triangles with rational area.

\begin{theorem}
If the quadrilateral is a three-sides-equal trapezoid $(a,a,c,a)$, then the torsion group is $\mathbb{Z}/6\mathbb{Z}$.
\end{theorem}

By scaling the sides of a rational cyclic quadrilateral, one can always assume that area is an integer $N$ (since only the quadrilaterals with rational area are considered). Using the definition of the generalized congruent number elliptic curve and taking into account Theorems 7.3, 7.4, and 7.5, one can restate Theorem 7.1 in the following way.
\begin{theorem}
Every non-square and non-three-sides equal trapezoidal rational cyclic quadrilateral with area $N \in \mathbb{N}$ gives rise to a generalized congruent number elliptic curve $E_{\alpha,-N^2}$ with positive rank.  Conversely, for any integer $N$ and generalized congruent curve $E_{\alpha,-N^2}$ with positive rank, there are infinitely many non-rectangular cyclic quadrilaterals.
\end{theorem}
Next the structure of torsion points are discussed.
\begin{pro}
There is no point $P$ on the curve $E_{\alpha,-n^2}$ such that $2P=(0,0)$. Consequently, the torsion group $T \neq \mathbb{Z}/4\mathbb{Z},$ $\mathbb{Z}/8\mathbb{Z}$, or $\mathbb{Z}/12\mathbb{Z}$.
\end{pro}

\begin{pro}
There are no points of order 5 on the curve $E_{\alpha,-n^2}$.
\end{pro}
By an old result, an elliptic curve with a rational point of order 5 will have its $j$-invariant of the form $(s^2+10s+5)^3/s$ for some $s \in \mathbb Q$ \cite{Fricke}.
Calculating the $j$-invariant of $E_{\alpha,-n^2}$, one must have $$256\frac{(\alpha^2+3n^2)^3}{n^4(\alpha^2+4n^2)}=\frac{(s^2+10s+5)^3}{s}.$$

\begin{pro}
The torsion group $T$ is not $\mathbb{Z}/2\mathbb{Z} \times \mathbb{Z}/8\mathbb{Z}$.
\end{pro}
\begin{pro}
The torsion group $T$ is not $\mathbb{Z}/2\mathbb{Z} \times \mathbb{Z}/6\mathbb{Z}$.
\end{pro}
Combining the above series of results, the following corollary is immediate.
\begin{cor}
Given a cyclic quadrilateral with corresponding elliptic curve $E_{\alpha,-n^2}$, the torsion group must be $\mathbb{Z}/2\mathbb{Z}$, $\mathbb{Z}/6\mathbb{Z}$,  $\mathbb{Z}/2\mathbb{Z} \times \mathbb{Z}/2\mathbb{Z}$, or $\mathbb{Z}/2\mathbb{Z} \times \mathbb{Z}/4\mathbb{Z}$.
\end{cor}

Note that all four torsion groups are possible.  As it was previously showed, any square will have torsion group $\mathbb{Z}/6\mathbb{Z}$.  For any $m>2$, and letting $a=m^2-4$ and $b=2m$, then the rectangle with side lengths $a$ and $b$ will have torsion group $\mathbb{Z}/2\mathbb{Z} \times \mathbb{Z}/2\mathbb{Z}$.

It turns out that if the cyclic quadrilateral is a non-square rectangle, two of the possible torsion groups can be ruled out.

\begin{pro}
The elliptic curve arising from a rectangle has a point of order 3 if and only if the rectangle is actually a square.
\end{pro}

\begin{cor}
Given a non-square rectangle with sides $a$ and $b$, then $T=\mathbb{Z}/2\mathbb{Z}$ or  $T=\mathbb{Z}/2\mathbb{Z} \times \mathbb{Z}/2\mathbb{Z}$. one gets the latter group if and only if $4a^2+b^2$ is a square.
\end{cor}

This characterization of congruent numbers can be added to the list of the many other known characterizations of congruent numbers. Several of these are given in Koblitz's book \cite{Koblitz}. Given an integer $n$, it is a well-known open problem to determine whether or not $n$ is congruent. A partial answer is given by Tunnell's theorem, which gives an easily testable criterion for determining if a number is congruent.  However, this result relies on the unproven Birch and Swinnerton-Dyer Conjecture for curves of the form $y^2=x^3-n^2x$.   The criterion involves counting the number of integral solutions $(x,y,z)$ to a few Diophantine equations of the form $ax^2+by^2+cz^2=n$.

In the remainder of the paper, the authors give two different new infinite families of congruent number elliptic curves with (at least) rank three.  Searching for families of congruent curves with high rank has been done before \cite{Duj5},\cite{Duj6},\cite{Rog},\cite{Rub}.  Currently, the best known results are a few infinite families with rank at least 3 \cite{Spe3},\cite{Rub}, and several individual curves with rank 7.

The authors conclude the paper by constructing an infinite family with rank at least 4 as well as an infinite family of rank 5. By computer search examples of the curves with rank 10 are also obtained.

\section{On Parametric Spaces of Bicentric Quadrilaterals}

In Euclidean geometry, a bicentric quadrilateral is a convex quadrilateral that has both a circumcircle passing through the four vertices and an incircle having the four sides as tangents.
In this work, the authors \cite{IKMZ} consider a bicentric quadrilateral with rational sides, and discuss the problem of finding such quadrilaterals where the ratio of the radii of the circumcircle and incircle is rational. The radii of these circles are denoted by $R$ and $r$ respectively. Examples of bicentric quadrilaterals  are squares, right kites, and isosceles tangential trapezoids.

First, we briefly recall some basic facts concerning these objects.
Let $ABCD$ be a convex quadrilateral with sides of rational length $a, b, c, d$ in anti-clockwise direction. It is bicentric if and only if the
opposite sides satisfy Pitot's theorem and the opposite angles are supplementary, that is, $s=a+c=b+d$, $\angle A+\angle C=\angle B+\angle D=\pi,$ where $s$ is the semi-perimeter.

In addition, we have
$$R=\frac{1}{4}\sqrt{\frac{(ab+cd)(ac+bd)(ad+bc)}{abcd}}$$
and $r=K/s$, where $K=\sqrt{abcd}$ is the area of the quadrilateral, see \cite[Page 67]{Gup}.

Consider a bicentric quadrilateral with rational sides. We discuss the problem of finding such quadrilaterals where the ratio of the radii of the circumcircle and incircle is rational. The authors show that this problem can be formulated in terms of a family of elliptic curves given by $E_a:y^2=x^3+(a^4-4a^3-2a^2-4a+1)x^2+16a^4x$ which have, in general,
\(\mathbb Z/8\mathbb Z\), and in rare cases \(\mathbb Z/2\mathbb Z\times\mathbb Z/8\mathbb Z\) as torsion subgroups. It can be show that the existence of infinitely many elliptic curves $E_a$ of rank at least two with torsion subgroup $\mathbb Z/8\mathbb Z$, parameterized by the points of an elliptic curve of rank at least one, and give five particular examples of rank $5$. It can be also show the existence of a subfamily of $E_a$ whose torsion subgroup is $\mathbb Z/2\mathbb Z\times\mathbb Z/8\mathbb Z$.

The authors show that a given bicentric quadrilateral, with $N$ rational, leads to a member of a family of elliptic curves, denoted by $E_a$, dependent on a single side.
They also discuss the torsion subgroup of $E_a(\mathbb Q)$, which is either $\mathbb Z/8\mathbb Z$ or $\mathbb Z/2\mathbb Z\times\mathbb Z/8\mathbb Z$.

The main results are the followings:

\begin{theorem}\label{ps}
All rational bicentric quadrilaterals with rational $N$ correspond to rational points on $E_a$ for some $a>0$, where $E_a(\mathbb Q)$ will have strictly positive rank.
\end{theorem}
\begin{theorem}\label{4.1}
The torsion subgroup of $E_a(\mathbb Q)$ equals $\mathbb Z/2\mathbb Z\times \mathbb Z/8\mathbb Z$ if and only if $$a=-\frac{r+1}{r(r-1)}\text{ for some }r\in\mathbb Q\setminus\{0,\pm1\}.$$
\end{theorem}
From this, refer  to \cite{GM,Du} to see what is known concerning the ranks of curves with prescribed torsion subgroup $\mathbb Z/2\mathbb Z\times \mathbb Z/8\mathbb Z$.

They close this section by listing, in Table \ref{T4}, 26 out of 27 known rank-three curves which come from the family $E_a$,  with $a$
defined as in Theorem \ref{4.1}, (see also \cite{Duj}).
\begin{table}[h]
\begin{center}
\begin{tabular}{lll}
\hline
$r$ & Authors & Date  \\ \hline
${12}/{17}$ & Connell-Dujella &  2000 \\
${47}/{18}$ &  Dujella &  2001  \\
${133}/{86}$ & Rathbun &  2003  \\
${201}/{239}$ & Campbell-Goins & 2003 \\
${299}/{589}$, $247/160$, $281/138$, & Dujella  &  2006  \\
$281/133$ && \\
$439/17$, $569/159$ & Rathbun & 2006 \\
$923/230$ & Dujella-Rathbun & 2006 \\
$247/419$, $200/99$, $337/65$ & Flores-Jones-Rollick-Weigandt-Rathbun & 2007 \\
$1017/352$ & Dujella & 2008 \\
${999}/{76}$, ${412}/{697}$, ${349}/{230}$,   & Fisher & 2009 \\
$217/425$, $440/217$, $309/470$,   &  & \\
$496/319$, $585/391$, $219/313$,  & & \\
$336/191$ && \\
${257}/{287}$ &  Rathbun  &  2013 \\ \hline
\end{tabular}
\vspace{.2cm}
\caption{Curves with torsion subgroup $\mathbb Z/2\mathbb Z\times\mathbb Z/8\mathbb Z$ having rank $3$}\label{T4}
\end{center}
\end{table}

\section{Elliptic Curves Arising from Brahmagupta Quadrilaterals}

A Brahmagupta quadrilateral is a cyclic quadrilateral whose sides, diagonals, and area are all integer values.
In this article, the authors \cite{IKMs} characterize the notions of Brahmagupta,
introduced  by  K. R. S. Sastry \cite{Sas}, by  means  of elliptic curves.
Motivated by these characterizations, The authors use Brahmagupta quadrilaterals to construct infinite families of elliptic curves with torsion group $\mathbb Z/2\mathbb Z\times \mathbb Z/2\mathbb Z$
having ranks (at least) 4, 5, and 6. Furthermore, by specializing they give examples from these families of specific curves with rank 9.

A cyclic polygon is one with vertices upon which a circle can be circumscribed.
Mathematicians have long been interested in Brahmagupta's work on Heron triangles and cyclic quadrilaterals. For example, consider Kummer's complex construction to generate
Heron quadrilaterals outlined in \cite{Dic}.
The existence and parametrization of quadrilaterals with rational side lengths (and additional conditions) has a long history \cite{A-N,Dai,Dic,Gup,I-V}.
Buchholz and Macdougall \cite{Buc-Mac1} have shown that there exist no nontrivial
cyclic quadrilaterals with rational area and having the property that the rational side lengths form an arithmetic or geometric progression.

Let us refer to a cyclic Heron quadrilateral as a Brahmagupta quadrilateral \cite{Sas}.
Sastry \cite{Sas} used Pythagorean triangles to construct general
Heron triangles and cyclic quadrilaterals whose side lengths, diagonals, and
area are integers, i.e., Brahmagupta quadrilaterals. He introduced a rational parametrization of the four sides of
these quadrilaterals:
\begin{equation}\label{1}
\left\{
\begin{array}{l}
a=(t(u+v)+1-uv)(u+v-t(1-uv)), \vspace{.2cm}\\
b=(1+u^2)(v-t)(1+tv), \vspace{.2cm}\\
c=t(1+u^2)(1+v^2), \vspace{.2cm}\\
d=(1+v^2)(u-t)(1+tu),
\end{array}
\right.
\end{equation}
where $t,u,v\in \mathbb Q$ such that $abcd \neq 0$.
Brahmagupta's formula gives the area $S$ of a cyclic quadrilateral, in terms of the side lengths $a, b, c,$ and $d$:
\begin{equation}
\label{area}
S=\sqrt{(s-a)(s-b)(s-c)(s-d)},
\end{equation}
where $s=(a+b+c+d)/2$.  Letting $d=0$, this reduces to the well known Heron's formula for the area of a triangle in terms of its side lengths.
Brahmagupta also determined formulas for the lengths of the diagonals:
\begin{equation}\label{diagonal}
D_1=\sqrt{\frac{(ac+bd)(ad+bc)}{ab+cd}}, \mbox{ and }
D_2=\sqrt{\frac{(ac+bd)(ab+cd)}{ad+bc}}.
\end{equation}
Using the parameterizations in \eqref{1} above, it is easily checked the area $S$ and diagonal lengths $D_1,D_2$ are rational.

A priori, there is no reason to associate Brahmagupta quadrilaterals with elliptic curves.  However, by the area formula (\ref{area}) one see the point $(\alpha,\beta)=(s,S)$ lies on the quartic
\begin{equation}
\label{area1}
\beta^2=(\alpha-a)(\alpha-b)(\alpha-c)(\alpha-d).
\end{equation}
This quartic is birationally equivalent to an elliptic curve in the following manner. Taking $\zeta=-1/\alpha$, the equation (\ref{area1}) turns into
$$
\beta^2=\left(a+\frac{1}{\zeta}\right)\left(b+\frac{1}{\zeta}\right)\left(c+\frac{1}{\zeta}\right)\left(d+\frac{1}{\zeta}\right),
$$
or equivalently
$$(\zeta^2\beta)^2=(1+\zeta a)(1+\zeta b)(1+\zeta c)(1+\zeta d).$$
By the substitution
$$x=\frac{(1+a\zeta)(d-b)(d-c)}{1+d\zeta},\quad y=\frac{\zeta^2\beta(d-a)(d-b)(d-c)}{(1+d\zeta)^2},$$
the curve  (\ref{area1}) thus turns into
\begin{equation}\label{2.2}
E:\quad y^2=x(x+(b-a)(d-c))(x+(c-a)(d-b)),
\end{equation}
or
\begin{equation}\label{2.3}
E:\quad y^2=x^3+Ax^2+Bx,
\end{equation}
where $A=(b-a)(d-c)+(c-a)(d-b)$ and $B=(b-a)(d-c)(c-a)(d-b)$.  Equation (\ref{2.2}) (or (\ref{2.3})) defines an elliptic curve so long as $a \neq b \neq c \neq d$.
Note that by setting $d=0$, this elliptic curve becomes the same elliptic curve studied in \cite{DP1}, which arose from Heron triangles. Since the coefficients of the the elliptic curve defined in terms of $a, b, c, d$ and these in turn depend on the variables $t, u, v$ by Sastry representation, it follows that the main family of elliptic curves defined over $\Bbb Q(t,u,v)$.

{Assuming $a \neq b \neq c \neq d$, the curve $E$ has three $2$-torsion points:
$$
 T_1=(0,0), \ T_2=((a-b)(d-c),0), \ T_3=((a-c)(d-b),0),
 $$
which shows the torsion group ${\mathcal T}$ contains $\Bbb Z/2\Bbb Z\times \Bbb Z/2\Bbb Z$.
It can be easily checked that ${\mathcal T}\simeq\Bbb Z/2\Bbb Z\times \Bbb Z/2\Bbb Z$ by using the specialization monomorphism \cite[III.11.4]{S2}.
Next, the authors construct an infinite family of rank at leat 3, and families of rank at least 4. From families of rank at least 4 they finally found a subfamily of rank at least 5 as ell as a family of rank at least 6.
A computer search also shows the following special examples for the ranks.

\begin{table}
\caption{Curves $E_{t,v,m}$ with high rank}\label{tab:1}
\centering
\begin{tabular}{  r  r  r || c}
$t$ & $v$ & $m$ & rank \\ \hline  \hline
2/5 & 3/4 & 3/2 & 9 \\ \hline
5/8 & $-3$ & 9/8 & 9 \\ \hline
2 & 1/8 & 7/6 & 9 \\ \hline
2 & 14 & 26 & 9 \\ \hline
14/5 & 10/11 & 13/5 & 9 \\  \hline
4 & $-3$ & 32 & 9 \\ \hline
4 & 2 & 87 & 9 \\ \hline
8 & 18 & 20 & $7\leq {\rm rank}\leq 9$ \\ \hline
10 & 1/4 & 7/4 & $8\leq {\rm rank}\leq 10$ \\ \hline
13 & $-3$ & 41 & 9 \\
\end{tabular}
\end{table}


\begin{thebibliography}{}
\bibitem{DP2} J. Aguirre, A. Dujella and J. C. Peral, On the rank of elliptic curves coming from rational Diophantine triples, Rocky Mountain J. of Math. 42, (2012), 1759-1776.
\bibitem{A-N} C. Alsina, and R. B. Nelsen, On the diagonals of a cyclic quadrilateral,
 Forum Geom., {\bf 7} (2007), 147--149.
\bibitem{Buc-Mac1} R.H. Buchholz, and J. A. Macdougall, Heron quadrilaterals with sides in sides in arithmetic or geometric progression, Bull. Aust. Math. Soc., {\bf59} (1999), 263--269.
\bibitem{GM} G. Campbell and E. H. Goins, Heron triangles, Diophantine problems and elliptic curves,
2003 (preprint).
\bibitem{Dai} L. Daia, On a conjecture, Gaz. Math., \textbf{89} (1984), 276-279.
\bibitem {Dic} L. E. Dickson, History of the theory of numbers II, Chelsea publishing company, New York, 1971.
\bibitem{Du} A. Dujella, On Mordell-Weil groups of elliptic curves induced by Diophantine triples,  Glas. Mat. Ser. III {\bf 42} (2007), 3--18.
\bibitem{Duj} A.  Dujella,   High rank  elliptic curves with prescribed torsion,\\
\url{https://web.math.pmf.unizg.hr/~duje/tors/tors.html}.
\bibitem{DP1} A. Dujella and J. C. Peral,  Elliptic curves coming from Heron triangles, Rocky Mountain J. Math. Vol. 44, No. 4 (2014), 1145-1160.
\bibitem{Duj6} A. Dujella,  A. S. Janfada, J. C. Peral, S. Salami, On the high rank $\pi/3$ and $2\pi/3$ - congruent number elliptic curves, Rocky Mountain J. Math. Vol. 44, No. 6, (2014), 1867-1880.
\bibitem{Duj5} A. Dujella, A. S. Janfada and S. Salami, A search for high rank congruent number elliptic curves, J. Integer Seq., \textbf{12} (2009), Article 09.5.8.
\bibitem{Fricke} R. Fricke, Die elliptischen Funktionen und ihre Anwendungen, Bd. 2. B. G. Teubner, Leipzig, (1922).
\bibitem{GM} E. H. Goin and D. Maddox, Heron triangles Via elliptic curves, Rocky Mountain J. of Math. Vol. 36. No. 5, 2006, 1511-1526.
\bibitem{Gup} R. C. Gupta, Parame\'{s}vara's rule for the circumradius of a cyclic quadrilateral, Historia Math., \textbf{4} (1977), no. 1, 67-74.
\bibitem{I-V} D. Ismailescu and A. Vojdany, Class preserving dissections of convex quadrilatrals, Forum Geom., {\bf 9} (2009), 195--211.
\bibitem{IKMZ}F. Izadi, F. Khoshnam, A. J. MacLeod, A. S. Zargar, On parametric spaces of bicentric quadrilaterals, Math. Slovaca 67 (2017), 611-622
\bibitem{IK} F. Izadi and F. Khoshnam, On elliptic curves Via Heron triangles and diophantine triples, J. of Mathematical Extension, Vol. 8, No. 3, (2014), 17-26.
\bibitem{IN} F. Izadi and  K. Nabadri, Diophantine equation $X^4+Y^4=2(U^4+V^4)$, Math. Slovaca 66 (2016), No.3 557-560.
\bibitem{IN2} F. Izadi and K. Nabardi, A family of elliptic curves with rank $\geq5$, Vol. 71, issue 2, 245-249.
\bibitem{IKM} F. Izadi, F. Khoshnam and D. Moody, Heron quadrilaterals Via elliptic curves, Rocky Mountain J. Math. Vol. 47, No. 4 (2017), 1227-1258.
\bibitem{IKMZ} F. Izadi, F. Khoshnam, D. Moody and A. S. Zargar, Elliptic curves arising from Brahmagupta quadrilaterals, Bull. Aust. Math. Soc, Vol. 92, (2015), 187-194.
\bibitem{IKZ1} F. Izadi, F. Khoshnam and A. S. Zargar, Rank of elliptic curves associated to Brahmagupta quadrilaterals, Colloquium Mathematicum, Vol. 143, No. 2, 2016, 187-192.
\bibitem{IKZ2} F. Izadi and A. S. Zargar,  A note on the Diophantine equations $x_1^k+x_2^k+x_3^k+x_4^k=2y_1^k+2y_2^k$, $k=3$, $6$, Vol. 20, 2014, No. 5, Pages 1—10.
\bibitem{Spe3} J. A. Johnstone and B. K. Spearman, Congruent number elliptic curves with
rank at least three, Canad. Math. Bull., \textbf{53} (2010), no. 4, 661-666.
\bibitem{Koblitz} N. Koblitz, Introduction to elliptic curves and modular forms, Springer-Verlag, New York, 1984.
\bibitem{Rog} N. F. Rogers, Rank computations for the congruent number elliptic curves. Exp. Math. \textbf{9}, (2000), no. 4, 591-594.
\bibitem {Rub} K. Rubin, and A. Silverberg, Rank frequencies for quadratic twists of elliptic curves, Exp. Math. \textbf{10} (2001), no. 4, 559-569.
\bibitem {Sas} K.R.S. Sastry, Brahmagupta quadrilaterals, Forum Geom., {\bf 2} (2002), 167--173.
\bibitem{S2} J. H. Silverman, Heights and the specialization map for families of abelian varieties, J. Reine Angew. Math. 342 (1983), 197-211.

\end{thebibliography}
\end{document}